\documentclass[11pt, a4paper]{article}

\usepackage{amsmath}
\usepackage{amssymb}
\usepackage{amsthm}
\usepackage[linktoc = all, hidelinks]{hyperref} 
\usepackage[nameinlink, noabbrev]{cleveref} 
\usepackage[margin = 2.5cm]{geometry}

\newcommand{\ext}[2]{\operatorname{ex}(#1, #2)}
\newcommand{\set}[1]{\mathopen{}\left\{#1\right\}\mathclose{}}

\usepackage{thmtools}
\usepackage{thm-restate}

\declaretheoremstyle[
spaceabove = .7\baselineskip\@plus.2\baselineskip\@minus.2\baselineskip, 
spacebelow = .7\baselineskip\@plus.2\baselineskip\@minus.2\baselineskip,
headfont = \normalfont\bfseries,
notefont = \mdseries, 
notebraces = {(}{)},
bodyfont = \normalfont\itshape,
postheadspace = .5em,
headpunct = .
]{bolditalic}

\declaretheoremstyle[
spaceabove = .5\baselineskip\@plus.2\baselineskip\@minus.2\baselineskip, 
spacebelow = .5\baselineskip\@plus.2\baselineskip\@minus.2\baselineskip,
headfont = \normalfont\itshape,
notefont = \mdseries, 
notebraces = {}{},
bodyfont = \normalfont,
postheadspace = .5em,
headpunct = .,
qed = \qedsymbol
]{proofstyle}

\declaretheorem[name = Theorem, numbered = yes, style = bolditalic, refname = {theorem,theorems}, Refname = {Theorem,Theorems}]{theorem}

\declaretheorem[name = Proof, numbered = no, style = proofstyle, refname = {proof,proofs}, Refname = {Proof,Proofs}]{Proof}

\title{\vspace{-0.7cm} A note on induced Tur\'{a}n numbers}
\author{Freddie Illingworth\thanks{DPMMS, University of Cambridge, UK. \emph{E-mail:} \textsf{fci21@cam.ac.uk}. Research supported by an EPSRC grant. \hfill \null \linebreak \emph{2020 MSC:} 05C35.}}
\date{}

\begin{document}
	
\maketitle
	
\begin{abstract}
	Loh, Tait, Timmons and Zhou introduced the notion of induced Tur\'{a}n numbers, defining $\ext{n}{\set{H, F\textnormal{-ind}}}$ to be the greatest number of edges in an $n$-vertex graph with no copy of $H$ and no induced copy of $F$. Their and subsequent work has focussed on $F$ being a complete bipartite graph. In this short note, we complement this focus by asymptotically determining the induced Tur\'{a}n number whenever $H$ is not bipartite and $F$ is not an independent set nor a complete bipartite graph. 
\end{abstract}

\section{Introduction}

Consider the following induced version of the classical Tur\'{a}n problem. What is the most edges in an $n$-vertex graph which does not contain $F$ as an induced subgraph? This question has a simple answer: if $F$ is a clique, then it is the original subgraph problem (and so answered by Tur\'{a}n's theorem~\cite{Turan1941}), and if $F$ is not complete, then the maximum is plainly $\binom{n}{2}$, as witnessed by the $n$-vertex complete graph $K_{n}$.

In an insightful paper, Loh, Tait, Timmons and Zhou~\cite{LTTZ2018} recovered an interesting and natural induced Tur\'{a}n problem by forbidding both an induced $F$ and a (not necessarily induced) copy of $H$. This removes the possibility of $K_{n}$ being extremal and so rules out the humdrum answer of $\binom{n}{2}$. To be precise, for a positive integer $n$ and graphs $H$ and $F$, they defined \emph{the induced Tur\'{a}n number}
\begin{equation*}
	\ext{n}{\set{H, F\textnormal{-ind}}},
\end{equation*}
to be the maximum number of edges in an $n$-vertex graph with no copy of $H$ and no induced copy of $F$.

In~\cite{LTTZ2018}, Loh, Tait, Timmons and Zhou focussed on complete bipartite $F$ proving general bounds for $\ext{n}{\set{H, K_{s, t}\textnormal{-ind}}}$ with some sharper results when $s = 2$. Nikiforov, Tait and Timmons~\cite{NTT2018} proved spectral improvements of their results. Ergemlidze, Gy\"{o}ri and Methuku~\cite{EGM2019} asymptotically determined $\ext{n}{\set{H, F\textnormal{-ind}}}$ in the important case where $H$ is an odd cycle and $F$ is $K_{2, t}$ or $K_{3, 3}$ (except when $H = C_{5}$ and $F = K_{2, 2}$ where they strengthened the upper bound). In~\cite{Illingworth2021K2t}, the author gave further improvements to the upper bounds for $\ext{n}{\set{H, K_{2, t}\textnormal{-ind}}}$.

All the work to date has focussed on complete bipartite $F$. The following theorem complements this focus.
\begin{theorem}\label{theorem}
	Let $H$ and $F$ be graphs with chromatic numbers $r + 1$ and $s + 1$ respectively. Then 
	\begin{equation*}
		\ext{n}{\set{H, F\textnormal{-ind}}} = \begin{cases}
			\bigl(1 - \tfrac{1}{r} + o(1)\bigr) \tbinom{n}{2} & \text{if } s \geqslant r \text{ or } F \text{ is not complete multipartite}, \\
			\bigl(1 - \tfrac{1}{s} + o(1)\bigr) \tbinom{n}{2} & \text{if } s < r \text{ and } F \text{ is complete multipartite}.
		\end{cases}
	\end{equation*}
\end{theorem}
This asymptotically determines the induced Tur\'{a}n number except if $H$ is bipartite or $F$ is an independent set or a complete bipartite graph.

\section{Proof of \texorpdfstring{\Cref{theorem}}{Theorem 1}}

We will use $T_{r}(n)$ to denoted the $r$-partite Tur\'{a}n graph on $n$ vertices. This is complete $r$-partite with parts as equal in size as possible. We remind the reader that $T_{r}(n)$ has $(1 - 1/r + o(1)) \binom{n}{2}$ edges.

\begin{Proof}[of \Cref{theorem}]
	If $s \geqslant r$ or $F$ is not complete multipartite, then $F$ is not an induced subgraph of the Tur\'{a}n graph $T_{r}(n)$. This also does not contain $H$, so
	\begin{equation*}
		\ext{n}{\set{H, F\textnormal{-ind}}} \geqslant e(T_{r}(n)) = \bigl(1 - \tfrac{1}{r} + o(1)\bigr) \tbinom{n}{2}.
	\end{equation*}
	On the other hand, the Erd\H{o}s-Simonovits theorem~\cite{ErdosSimonovits1966} says that any $n$-vertex graph not containing $H$ has at most $(1 - 1/r + o(1)) \binom{n}{2}$ edges. This completes the first part of the theorem.
	
	Now suppose that $s < r$ and $F$ is complete $(s + 1)$-partite. Firstly the Tur\'{a}n graph $T_{s}(n)$ contains neither $H$ nor $F$ as subgraphs (let alone induced ones) so,
	\begin{equation*}
		\ext{n}{\set{H, F\textnormal{-ind}}} \geqslant e(T_{s}(n)) = \bigl(1 - \tfrac{1}{s} + o(1)\bigr) \tbinom{n}{2}.
	\end{equation*}
	As $F$ is complete $(s + 1)$-partite, there is a positive integer $t$ such that $F$ is an induced subgraph of $K_{s + 1}(t) = T_{s + 1}(t(s + 1))$. By Ramsey's theorem~\cite{Ramsey1930}, there is a positive integer $m$ such that any $m$-vertex graph contains either an independent $t$-set or a copy of $H$. 
	
	Fix $\varepsilon > 0$ and let $G$ be an $n$-vertex graph with at least $(1 - 1/s + \varepsilon) \binom{n}{2}$ edges where $n$ is sufficiently large. By the Erd\H{o}s-Stone theorem~\cite{ErdosStone1946}, $G$ contains a copy of $K_{s + 1}(m)$. Let the parts of the $K_{s + 1}(m)$ be $V_{1}, \dotsc, V_{s + 1}$ so each one has $m$ vertices. If any $V_{i}$ contains $H$, then $G$ does. Otherwise, by the definition of $m$, each $V_{i}$ contains an independent set of size $t$. Thus $G$ contains an induced copy of $K_{s + 1}(t)$ and so an induced copy of $F$. Thus, for all large $n$,
	\begin{equation*}
		\ext{n}{\set{H, F\textnormal{-ind}}} \leqslant \bigl(1 - \tfrac{1}{s} + \varepsilon\bigr) \tbinom{n}{2}.
	\end{equation*}
	This holds for arbitrary $\varepsilon > 0$, as required.
\end{Proof}

\end{document}